\documentclass[tbtags,reqno]{amsart}
\usepackage{epsfig,amsthm,amsopn,amsfonts,amssymb,bbm}

\newtheorem{theorem}{Theorem}
\newtheorem{lemma}[theorem]{Lemma}
\newtheorem{proposition}[theorem]{Proposition}
\newtheorem{example}[theorem]{Example}

\newtheorem{corollary}[theorem]{Corollary}

\newtheorem{definition}[theorem]{Definition}

\newtheorem{property}[theorem]{Property}

\theoremstyle{remark}
\newtheorem{remark}[theorem]{Remark}

\def \la {\lambda}

\def\om {\omega}
\def\kbnd {\mathfrak c^{-1}}
\def\core {\mathfrak c}

\def \core {{\mathfrak c}}
\def \rcore {{\mathfrak r}}
\def\la {\lambda}
\def\gg {\gamma}
\def\aa {\alpha}

\edef\savecatcodeat{\the\catcode`@}
\catcode`\@=11

\def\tb@ifSpecChars#1#2{#1}
\def\tb@ifNoSpecChars#1#2{#2}

\def\tableau{%
  \bgroup
  \@ifstar{\let\Tif\tb@ifNoSpecChars\tb@tableauB}
          {\let\Tif\tb@ifSpecChars\tb@tableauB}}

\def\tb@tableauB{
  \@ifnextchar[{\tb@tableauC}{\tb@tableauC[]}}

\def\tb@tableauC[#1]{\hbox\bgroup%
    \let\\=\cr
    \def\bl{\global\let\tbcellF\tb@cellNF}%
    \def\tf{\global\let\tbcellF\tb@cellH}
    \dimen2=\ht\strutbox \advance\dimen2 by\dp\strutbox%
    \ifx\baselinestretch\undefined\relax%
    \else%
       \dimen0=100sp \dimen0=\baselinestretch\dimen0%
       \dimen2=100\dimen2 \divide\dimen2 by\dimen0%
    \fi%
    \let\tpos\tb@vcenter
    \tb@initYoung
    \tb@options#1\eoo
    \let\arrow\tb@arrow%
    \dimen0=\Tscale\dimen2%
    \dimen1=\dimen0 \advance\dimen1 by \tb@fframe%
    \lineskip=0pt\baselineskip=0pt
    \def\tb@nothing{}%
    \def\endcellno{$\rss\egroup\bss\egroup}
    \def\endcell{\endcellno\kern-\dimen0}
    \def\begincell{\vbox to\dimen0\bgroup\vss\hbox to\dimen0\bgroup\hss$}%
    \let\overlay\tb@overlay%
    \let\fl\tb@fl%
    \let\lss\hss\let\rss\hss\let\tss\vss\let\bss\vss
    \def\mkcell##1{
        \let\tbcellF\tb@cellD
        \def\tb@cellarg{##1}
        \ifx\tb@cellarg\tb@nothing\let\tb@cellarg\tb@cellE\fi%
        \begincell\tb@cellarg\endcellno
        \tbcellF}
    \let\savecellF\tbcellF
     \Tif{\catcode`,=4\catcode`|=\active}{}\tb@tableauD}%

\let\tb@savetableauD\tableauD
{
    \catcode`|=\active \catcode`*=\active \catcode`~=\active%
    \catcode`@=\active
\gdef\tableauD#1{%
  \Tif{
    \mathcode`|="8000 \mathcode`*="8000%
    \mathcode`~="8000 \mathcode`@="8000%
    \def@{\bullet}%
    \let|\cr
    \let*\tf
    \let~\sk
  }{}%
  \tpos{\tabskip=0pt\halign{&\mkcell{##}\cr#1\crcr}}%
  \global\let\tbcellF\savecellF
  \egroup
  \egroup}
}
\let\tb@tableauD\tableauD
\let\tableauD\tb@savetableauD
\let\tb@savetableauD\undefined

\def\tb@options#1{\ifx#1\eoo\relax\else\tb@option#1\expandafter\tb@options\fi}

\def\tb@option#1{%
  \if#1t\let\tpos\tb@vtop\fi
  \if#1c\let\tpos\tb@vcenter\fi
  \if#1b\let\tpos\vbox\fi
  \if#1F\tb@initFerrers\fi
  \if#1Y\tb@initYoung\fi
  \if#1s\tb@initSmall\fi
  \if#1m\tb@initMedium\fi
  \if#1l\tb@initLarge\fi
  \if#1p\tb@initPartition\fi
  \if#1a\tb@initArrow\fi
}

\def\tb@vcenter#1{\ifmmode\vcenter{#1}\else$\vcenter{#1}$\fi}

\def\tb@vtop#1{\hbox{\raise\ht\strutbox\hbox{\lower\dimen0\vtop{#1}}}}

\def\tb@initPartition{\def\Tscale{.3}}
\def\tb@initSmall{\def\Tscale{1}}
\def\tb@initMedium{\def\Tscale{2}}
\def\tb@initLarge{\def\Tscale{3}}

\def\tb@initArrow{\dimen2=1.25em}

\def\tb@initYoung{%
  \def\tb@cellE{}
  \let\tb@cellD\tb@cellN
  \def\sk{\global\let\tbcellF\tb@cellNF}}
\def\tb@initFerrers{%
  \def\tb@cellE{\bullet}
  \let\tb@cellD\tb@cellNF
  \def\sk{\bullet}}

\tb@initMedium

\def\tb@sframe#1{%
  \vbox to0pt{
    \vss
    \hbox to0pt{%
      \hss
      \vbox to\dimen1{
        \hrule depth #1 height0pt
        \vss
        \hbox to\dimen1{
          \vrule width #1 height\dimen1
          \hss
          \vrule width #1
          }%
        \vss
        \hrule height #1 depth 0in
        }%
      \kern-\tb@hframe
      }%
    \kern-\tb@hframe}}

\def\tb@hframe{.2pt}\def\tb@fframe{.4pt}\def\tb@bframe{2pt}
\def\tb@cellH{\tb@sframe{\tb@bframe}}       
\def\tb@cellNF{}                            
\def\tb@cellN{\tb@sframe{\tb@fframe}}       
\let\tbcellF\tb@cellN                       

\def\tb@rpad{1pt}
\def\tb@lpad{1pt}
\def\tb@tpad{1.8pt}
\def\tb@bpad{1.8pt}

\def\tb@overlay{\endcell\@ifnextchar[{\tb@overlaya}{\begincell}}
\def\tb@overlaya[#1]{\vbox to\dimen0\bgroup%
  \tb@overlayoptions#1\eoo%
  \tss\hbox to\dimen0\bgroup\lss$}
\def\tb@overlayoptions#1{\ifx#1\eoo\relax\else\tb@overlayoption#1\expandafter\tb@overlayoptions\fi}

\def\tb@overlayoption#1{
  \if#1t\def\tss{\vskip\tb@tpad}\let\bss\vss\fi
  \if#1c\let\tss\vss\let\bss\vss\fi
  \if#1b\def\bss{\vskip\tb@bpad}\let\tss\vss\fi
  \if#1l\def\lss{\hskip\tb@lpad}\let\rss\hss\fi
  \if#1m\let\lss\hss\let\rss\hss\fi
  \if#1r\def\rss{\hskip\tb@rpad}\let\lss\hss\fi
}

\def\tb@fl{\endcell\begincell\vrule depth 0pt width \dimen0 height \dimen0 \endcell\begincell}

\@ifundefined{diagram}{}{
\def\tb@arrowpad{.5}

\newoptcommand{\tb@arrow}{\@ne}[2]{%
  \endcell
   \begingroup%
   \let\dg@getnodesize\tb@getnodesize
   \dg@USERSIZE=#1\relax%
   \ifnum\dg@USERSIZE<\@ne \dg@USERSIZE=\@ne \fi%
   \dg@parse{#2}%
   \dg@label{\tb@draw{#1}{#2}}}

\def\tb@getnodesize#1#2#3#4#5{\dimen3=\tb@arrowpad\dimen2 #4=\dimen3 #5=\dimen3\relax}
\def\tb@getnodesize#1#2#3#4#5{\ifnum#2=0\ifnum#3=0\tb@getnodesizetail{#4}{#5}\else\tb@getnodesizehead{#4}{#5}\fi\else\tb@getnodesizehead{#4}{#5}\fi}
\def\tb@getnodesizetail#1#2{\dimen3=.5\dimen2 #1=\dimen3 #2=\dimen3}
\def\tb@getnodesizehead#1#2{\dimen3=.5\dimen2 #1=\dimen3 #2=\dimen3}

\def\tb@draw#1#2#3#4{%
        \dg@X=0\dg@Y=0\dg@XGRID=1\dg@YGRID=1\unitlength=.001\dimen0%
        \dg@LBLOFF=\dgLABELOFFSET \divide\dg@LBLOFF\unitlength%
        \dg@drawcalc
        \begincell
        \let\lams@arrow\tb@lams@arrow
        \begin{picture}(0,0)\begingroup\dg@draw{#1}{#2}{#3}{#4}\end{picture}%
        \endcell
        \endgroup
        \begincell}
}

\def\tb@lams@arrow#1#2{%
 \lams@firstx\z@\lams@firsty\z@
 \lams@lastx#1\relax\lams@lasty#2\relax
 \lams@center\z@
 \N@false\E@false\H@false\V@false
 \ifdim\lams@lastx>\z@\E@true\fi
 \ifdim\lams@lastx=\z@\V@true\fi
 \ifdim\lams@lasty>\z@\N@true\fi
 \ifdim\lams@lasty=\z@\H@true\fi
 \NESW@false
 \ifN@\ifE@\NESW@true\fi\else\ifE@\else\NESW@true\fi\fi
 \ifH@\else\ifV@\else
  \lams@slope
  \ifnum\lams@tani>\lams@tanii
   \lams@ht\ten@\p@\lams@wd\ten@\p@
   \multiply\lams@wd\lams@tanii\divide\lams@wd\lams@tani
  \else
   \lams@wd\ten@\p@\lams@ht\ten@\p@
   \divide\lams@ht\lams@tanii\multiply\lams@ht\lams@tani
  \fi
 \fi\fi
 \ifH@  \lams@harrow
 \else\ifV@ \lams@varrow
 \else \lams@darrow
 \fi\fi
}

\catcode`\@=\savecatcodeat
\let\savecatcodeat\undefined

\begin{document}

\title[Quantum cohomology and the $k$-Schur basis]
{Quantum cohomology and the $k$-Schur basis}

\author{Luc Lapointe}
\thanks{Research supported in part by FONDECYT (Chile) grant \#1030114,
the Programa Formas Cuadr\'aticas of the Universidad de Talca,
and NSERC (Canada) grant \#250904}
\address{Instituto de Matem\'atica y F\'{\i}sica,
Universidad de Talca, Casilla 747, Talca, Chile}
\email{lapointe@inst-mat.utalca.cl}

\author{Jennifer Morse}
\thanks{Research supported in part by NSF grant \#DMS-0400628}
\address{Department of Mathematics,
University of Miami, Coral Gables, Fl 33124}
\email{morsej@math.miami.edu}


\begin{abstract}
We prove that structure constants related to Hecke algebras at 
roots of unity are special cases of $k$-Littlewood-Richardson coefficients 
associated to a product of $k$-Schur functions.  As a consequence, both the 
3-point Gromov-Witten invariants appearing in the quantum cohomology of 
the Grassmannian, and the fusion coefficients for the WZW conformal field theories
associated to $\widehat{su}(\ell)$ are shown to be $k$-Littlewood Richardson 
coefficients.  From this, Mark Shimozono conjectured that the $k$-Schur functions 
form the Schubert basis for the homology of the loop Grassmannian, whereas
$k$-Schur coproducts correspond to the integral cohomology of the loop 
Grassmannian.  We introduce dual $k$-Schur functions defined on weights 
of $k$-tableaux that, given Shimozono's conjecture, form the Schubert basis 
for the cohomology of the loop Grassmannian.  We derive several properties of
these functions that extend those of skew Schur functions. 
\end{abstract}

\maketitle


\section{Introduction}

The study of Macdonald polynomials led to the discovery of 
symmetric functions, $s_\lambda^{(k)}$, indexed by partitions
whose first part is no larger than a fixed integer $k\geq 1$.
Experimentation suggested that these functions play the 
fundamental combinatorial role of the Schur basis in 
the symmetric function subspace $\Lambda^{k}=\mathbb Z[h_1,\ldots,h_k]$;
that is, they satisfy properties generalizing 
classical properties of Schur functions such 
as Pieri and Littlewood-Richardson rules.  The study of the 
$s_\lambda^{(k)}$ led to several 
different characterizations \cite{[LLM],[LMfil],[LMkschur]} (conjecturally 
equivalent) and to the proof of many of these combinatorial conjectures.  
We thus generically call the functions $k$-Schur functions, but in
this article consider only the definition presented in \cite{[LMkschur]}.

\medskip

Although prior work with $k$-Schur functions concentrated on proving that 
they act as the ``Schur basis" for $\Lambda^k$, the analogy was so striking 
that it seemed likely to extend beyond combinatorics to fields such as 
algebraic geometry and representation theory.  Our main finding in this 
direction is that the $k$-Schur functions are connected to representations 
of Hecke algebras $H_\infty(q)$, where $q$ is a root of unity, and they
provide the natural basis for work in the quantum cohomology of the 
Grassmannian just as the Schur functions do for the usual cohomology.
In particular, the 3-point Gromov-Witten invariants are none other than 
relevant cases of ``$k$-Littlewood-Richardson coefficients",
the expansion coefficients in
\begin{equation}
s_\lambda^{(k)} s_\mu^{(k)} =\sum_{\nu: \nu_1\leq k}
c_{\lambda\mu}^{\nu , k} s_\nu^{(k)}\,.
\end{equation}

\medskip

To be precise, in Schubert calculus, the cohomology ring of the Grassmannian 
$Gr_{\ell n}$ (the manifold of $\ell$-dimensional subspaces of $\mathbb C^n$)
has a basis given by Schubert classes $\sigma_\lambda$ that are
indexed by partitions $\lambda\in \mathcal P^{\ell n}$ 
that fit inside an $\ell \times (n-\ell)$ rectangle.
There is an isomorphism,
$$
H^*(Gr_{\ell n})\cong \Lambda^{\ell}/\langle e_{n-\ell+1},\ldots,e_n\rangle\, ,
$$
where the Schur function $s_\lambda$ maps to the Schubert class
$\sigma_\lambda$ when $\lambda\in\mathcal P^{\ell n}$.  Since 
$s_{\lambda}$ is zero modulo the ideal when 
$\lambda \not \in \mathcal P^{\ell n}$, the structure constants 
of $H^*(Gr_{\ell n})$ in the basis of Schubert classes:
$$
\sigma_\lambda\sigma_\mu = \sum_{\nu\in\mathcal P^{\ell n}}
c_{\lambda\mu}^\nu \sigma_\nu
\,,
$$
can be obtained from the Littlewood-Richardson coefficients for 
Schur functions,
$$
s_\lambda s_\mu = 
\sum_{\nu\in\mathcal P^{\ell n}}
c_{\lambda\mu}^\nu s_\nu
+\sum_{\nu\not\in\mathcal P^{\ell n}}
c_{\lambda\mu}^\nu s_\nu
\, ,
$$
which have well known combinatorial interpretations.

\medskip

The small quantum cohomology ring of the Grassmannian $QH^*(Gr_{\ell n})$
is a deformation of the usual cohomology that has become the object of 
much recent attention (e.g. \cite{[Ag],[Wi]}).  As a linear space, this is 
the tensor product 
$H^*(Gr_{\ell n})\otimes\mathbb Z[q]$ and the $\sigma_\lambda$ with
$\lambda\in \mathcal P^{\ell n}$ form a $\mathbb Z[q]$-linear basis 
of $QH^*(Gr_{\ell n})$.  Multiplication is a $q$-deformation of the product
in $H^*(Gr_{\ell n})$, defined by
$$
\sigma_\lambda * \sigma_\mu = \sum_{\nu\in\mathcal P^{\ell n}\atop
|\nu|=|\lambda|+|\mu|-dn}
q^dC_{\lambda\mu}^{\nu,d}\sigma_\nu
\,.
$$
The $C_{\lambda\mu}^{\nu,d}$ are the 3-point Gromov-Witten invariants,
which count the number of rational curves of degree $d$ in $Gr_{\ell n}$
that meet generic translates of the Schubert varieties associated
to $\lambda,\mu$, and $\nu$.  Finding a combinatorial interpretation for 
these constants is an interesting open problem that would have applications 
to many areas, including the study
of the Verlinde fusion algebra \cite{[TUY]} as well as the computation
of certain knot invariants \cite{[Tur]}.

\medskip

As with the usual cohomology, quantum cohomology can be connected to
symmetric functions by:
$$
QH^*(Gr_{\ell n}) \cong \left(\Lambda^{\ell}\otimes\mathbb Z[q]\right)/
J^{\ell n}_q \, ,
$$
where $J^{\ell n}_q=\langle e_{n-\ell+1},\ldots,e_{n-1},e_n+(-1)^{\ell}q\rangle$.
When $\lambda \in \mathcal P^{\ell n}$, the Schubert class $\sigma_\lambda$
still maps to the Schur function $s_\lambda$, but unfortunately 
when $\lambda \not \in \mathcal P^{\ell n}$, some $s_{\lambda}$ 
are not zero modulo the ideal.  Thus, the Schur functions cannot be used to
directly obtain the quantum structure constants. Instead, these Gromov-Witten
invariants arise as the expansion coefficients in
$$
s_\lambda \, s_\mu = 
\sum_{\nu\in\mathcal P^{\ell n}\atop
|\nu|=|\lambda|+|\mu|-dn}
q^d C_{\lambda\mu}^{\nu,d}s_\nu 
\mod J^{\ell n}_q
\, ,
$$
and to compute the coefficients, an algorithm involving negatives 
\cite{[Cu],[Kac],[Wa]} must be used to reduce a Schur function  
modulo the ideal $J^{\ell n}_q$.

\medskip

Remarkably, by first working with an ideal that arises in the context of 
Hecke algebras at roots of unity, we find that the $k$-Schur functions 
circumvent this problem:  a $k$-Schur function maps to a single Schur 
function times a $q$ power 
(with no negatives) or to zero, modulo the ideal.  To be more precise, 
let $I^{\ell n}$ denote the ideal
$$
I^{\ell n}= \Bigl \langle s_{\lambda} \, 
\Big| \, 
 \# \{ j \, |  \, \lambda_j < \ell\}=n-\ell+1 \Bigr \rangle \, .
$$
A basis for $\Lambda^{\ell}/I^{\ell n}$ is given by the  
Schur functions indexed by partitions in $\Pi^{\ell n}$, 
the set of partitions with no part larger than $\ell$ and 
no more than $n-\ell$ rows of length smaller than $\ell$.
In \cite{[GW]}, certain structure constants
associated to representations of Hecke algebras at roots of unity are shown to 
be the expansion coefficients in
$$
s_\lambda \, s_\mu = \sum_{\nu\in \Pi^{\ell n}}
a_{\lambda\mu}^{\nu} s_{\nu} \mod I^{\ell n}\, .
$$
We prove that the $a_{\lambda\mu}^{\nu}$ are just special cases of 
$k$-Littlewood-Richardson coefficients by showing that when 
$\nu\in \Pi^{\ell n}$, the $k$-Schur function $s_{\nu}^{(k=n-1)}$ 
modulo the ideal $I^{\ell n }$ is simply $s_{\nu}$, and is zero 
otherwise.  Thus it is revealed that the  $a_{\lambda \mu}^{\nu}$
are coefficients in the expansion:
$$
s_\lambda^{(k)} \, s_\mu^{(k)} = \sum_{\nu\in \Pi^{\ell n}}
a_{\lambda\mu}^{\nu} \, s_{\nu}^{(k)} +
\sum_{\nu\not\in \Pi^{\ell n}}
c_{\lambda \mu}^{\nu,k} s_{\nu}^{(k)}
\,.
$$

\medskip

We can then obtain the 3-point Gromov-Witten invariants from this result by simply
computing $s_{\nu}$ modulo $J_{\ell n}^q$ for $\nu \in \Pi^{\ell n}$, since
$I^{\ell n}$ is a subideal of $J^q_{\ell n}$.  In this case, $s_{\nu}$ 
beautifully reduces to positive $s_{\mathfrak r (\nu)}$ times a $q$ power,
where $\mathfrak r (\nu)$ is the $n$-core of $\nu$.  Consequently, we
prove that the 3-point Gromov-Witten invariants are none other than 
certain $k$-Schur function Littlewood Richardson coefficients.  To
be more specific,
$$
C_{\lambda \mu}^{\nu,d} = c_{\lambda \mu}^{\hat \nu, n-1} \, ,
$$
where the value of $d$ associates a unique element $\hat \nu\in\Pi^{\ell n}$ 
(given explicitly in Theorem~\ref{gromovwitten}) to each 
$\nu\in\mathcal P^{\ell n}$. 

\medskip

It also follows from our results that the $k$-Littlewood-Richardson 
coefficients include the fusion rules for the Wess-Zumino-Witten 
conformal field theories associated to $\widehat{su}(\ell)$ at level $n-\ell$, 
since the algorithm given by Kac \cite{[Kac]} and Walton \cite{[Wa]} for 
computing in the fusion algebra reduces to the one given by Goodman and 
Wenzl \cite{[GW]} for computing the Hecke algebra structure constants.

\medskip

It is important to note that since the Gromov-Witten invariants under
consideration are indexed by partitions fitting inside a rectangle, they are given 
by only a subset of the $k$-Littlewood-Richardson coefficients. We thus naturally 
sought the larger picture that would be explained by the complete set of 
$k$-Littlewood Richardson coefficients.  In discussion with Mark Shimozono 
about this problem, he conjectured that the $k$-Schur functions form the 
Schubert basis for the homology of the affine (loop) Grassmannian of $GL_{k+1}$, and
that the $k$-Schur expansion coefficients of the $k$-Schur coproduct give the 
integral cohomology of the loop Grassmannian.  Here we introduce a family of 
functions dual to the $k$-Schur functions, defined by the
weight of certain ``$k$-tableaux" related to the affine symmetric
group \cite{[LMcore]}.  Following the theory of skew Schur functions,
we prove a number of results about these dual $k$-Schur functions
including that their symmetry relies on a generalization \cite{[LMkschur]} 
of the Bender-Knuth involution \cite{[BK]}.
In particular, we show that the coefficients in a product of dual $k$-Schur 
functions are the structure constants in the $k$-Schur coproduct, 
implying from Shimozono's conjecture that the dual $k$-Schur functions 
form the Schubert basis for the cohomology of the loop Grassmannian.  

\medskip

In addition to using the $k$-Schur functions to study the Gromov-Witten invariants 
and the loop Grassmannian, they are a natural tool to seek ``affine Schubert 
polynomials".  Our results strongly support the idea of Michelle Wachs that 
the (dual) $k$-Schur functions provide the symmetric Grassmannian component 
of a larger family of polynomials that are analogous to Schubert polynomials,
but indexed instead by affine permutations.  After discussion with Thomas Lam 
of the work presented here, he made a beautiful step in this direction by
introducing a family of ``affine Stanley symmetric 
functions" that reduce in special cases to the dual $k$-Schur functions (called ``affine 
Schur functions" in \cite{[Lam]}).  Details of a connection between the 
dual $k$-Schur functions and the cylindric Schur functions of 
\cite{[Po]} is also carried out in \cite{[Lam]}.

\section{definitions}

Let $\Lambda$ denote the ring of symmetric functions, generated by the 
elementary symmetric functions $e_r=\sum_{i_1<\ldots <i_r}x_{i_1}
\cdots x_{i_r}$, or equivalently by the complete symmetric functions
$h_r=\sum_{i_1\leq\ldots\leq i_r}x_{i_1}\cdots x_{i_r}$, and let
$\Lambda^{k} = \mathbb Z[h_1,\ldots,h_k]$.  Bases for $\Lambda$ are indexed 
by partitions $\lambda=(\lambda_1\geq\dots\geq\lambda_m>0)$ whose 
degree $\lambda$ is $|\lambda|=\lambda_1 +\cdots +\lambda_m$
and whose length $\ell(\lambda)$ is the number of parts $m$.
Each partition $\lambda$ has an associated Ferrers diagram
with $\lambda_i$ lattice squares in the $i^{th}$ row,
from the bottom to top. Any lattice square in the Ferrers diagram 
is called a cell, where the cell $(i,j)$ is in the $i$th row and $j$th 
column of the diagram.  Given  a partition $\lambda$, its 
conjugate $\lambda'$ is the diagram 
obtained by reflecting  $\lambda$ about the main diagonal.  A partition 
$\lambda$ is ``{\it $k$-bounded}'' if $\lambda_1 \leq k$ and the set
of all such partitions is denoted $\mathcal P^k$.  The set 
$\mathcal P^{\ell n}$ is the partitions fitting inside an 
$\ell \times (n-\ell)$ rectangle (with $n-\ell$ rows of size $\ell$).  
We say that $\lambda \subseteq \mu$ when $\lambda_i \leq \mu_i$ for
all $i$.  Dominance order $\unrhd$ on partitions is defined by 
$\lambda\unrhd\mu$ when $\lambda_1+\cdots+\lambda_i\geq
\mu_1+\cdots+\mu_i$ for all $i$, and $|\lambda|=|\mu|$.

\medskip

More generally, for $\rho \subseteq \gamma$,
the skew shape $\gg/\rho$ is identified with its
diagram $\{(i,j) : \rho_i<j\leq \gg_i\}$.
Lattice squares that do not lie in $\gg/\rho$ 
will be simply called ``{\it squares$\, $}''. We say that any 
$c\in \rho$ lies ``{\it below $\, $}'' $\gg/\rho$.
The ``{\it hook $\, $}'' of any lattice square $s\in \gg$ is defined as the 
collection of cells  of $\gg/\rho$ that lie inside the $L$ with $s$ as 
its corner. This is intended to  apply to all  $s\in \gg$ including those 
below $\gg/\rho$.  For example, the hook of $s=(1,3)$ is  
depicted by the framed cells:
\begin{equation}
\gg/\rho\,=\,(5,5,4,1)/(4,2)\,=\,
{\tiny{\tableau*[scY]
{\cr&&\tf&\cr\bl&\bl&\tf&&\cr\bl &\bl &\bl s &\bl & \tf}}} \, .
\end{equation}
The {\it ``hook-length"} of $s$, $h_s(\gg/\rho)$, is the number of 
cells in the hook of $s$.
In the preceding example, $h_{(1,3)}\big((5,5,4,1)/(4,2)\big)=3$ 
and $h_{(3,2)}\big((5,5,4,1)/(4,2)\big)=3$.  
A cell or square has a $k$-bounded hook if it's hook-length is
no larger than $k$.

\medskip

A {\it ``$p$-core''} is a partition that does not contain any hooks of 
length $p$, and $\mathcal C^p$ will denote the set of all $p$-cores.
The ``{\it $p$-residue}'' of square  $(i,j)$ is $j-i \mod p$;
that is, the label of this square when squares are periodically labeled 
with $0,1,\ldots,p-1$, where zeros lie on the main diagonal
(see \cite{[JK]} for more on cores and residues).  The 5-residues 
associated to the 5-core $(6,4,3,1,1,1)$ are
$$
{\tiny{\tableau[scY]{\bl 4 |0|1|2,\bl 3|3,4,0,\bl 1|4,0,1,2,
\bl 3| 0,1,2,3,4,0,\bl 1}}}
$$

\medskip

A {\it ``tableau"} is a filling of a Ferrers shape with integers 
that strictly increase in columns and weakly increase in rows.  
The {\it ``weight"} of a given tableau is the composition $\alpha$ where 
$\alpha_i$ is the multiplicity of $i$ in the tableau.
A ``{\it Schur function}" can be defined by
\begin{equation}
s_{\lambda}= \sum_{T} x^{T} \, ,
\end{equation}
where the sum is over all tableaux of shape $\lambda$,
and where $x^T=x^{\text{weight}(T)}$.

\def \CC {{\mathcal C}}
\def \CP {{\mathcal P}}
\def \CCk {\CC_{k+1}}

\section{$k$-Schur functions} \label{sect3}

There are several conjecturally equivalent characterizations for the
$k$-Schur functions.  Here we use the definition explored in \cite{[LMkschur]}
that relies on a family of tableaux related to the affine symmetric group.

\begin{definition}
\cite{[LMcore]}
\label{defktabgen} 
Let $\gg$ be a $k+1$-core, $m$ be the number of $k$-bounded hooks 
of $\gg$, and $\aa=(\aa_1,\ldots,\aa_r)$ be a composition of
$m$. A ``$k$-tableau" of shape $\gg$ and ``$k$-weight" $\aa$ 
is a filling of $\gg$ with integers $1,2,\ldots,r$ such that

\noindent
(i) rows are weakly increasing and columns are strictly increasing

\noindent
(ii) the collection of cells filled with letter $i$ are labeled with exactly 
$\alpha_i$ distinct $k+1$-residues.
\end{definition}

\begin{example}
\label{exssktab}
The $3$-tableaux of $3$-weight $(1,3,1,2,1,1)$ and shape $(8,5,2,1)$ are:
\begin{equation}
{\tiny{\tableau*[scY]{5\cr 4&6\cr2&3&4&4&6\cr 1&2&2&2&3&4&4&6 }}} \quad
{\tiny{\tableau*[scY]{6\cr 4&5\cr 2&3&4&4&5\cr 1&2&2&2&3&4&4&5 }}} \quad
{\tiny{\tableau*[scY]{4\cr 3&6\cr 2&4&4&5&6\cr 1&2&2&2&4&4&5&6 }}}
\end{equation}
\end{example}

\medskip

The definition of $k$-tableaux easily extends.

\begin{definition}
\label{defktabgenskew}
Let $\delta\subseteq\gg$ be $k+1$-cores with $m_1$ and $m_2$
$k$-bounded hooks respectively, and  let
$\aa=(\aa_1,\ldots,\aa_r)$ be a composition of $m_1-m_2$.  A
{\it ``skew $k$-tableau"} of shape $\gg/\delta$ and ``$k$-weight" $\aa$
is a filling of $\gg/\delta$ with integers $1,2,\ldots,r$ such that

\smallskip
\noindent
(i) rows are weakly increasing and columns are strictly increasing

\smallskip
\noindent
(ii) the collection of cells filled with letter $i$ are labeled by exactly
$\alpha_i$ distinct $k+1$-residues.
\end{definition}

\medskip

Although a $k$-tableau is associated to a shape $\gamma$ and weight
$\alpha$, in contrast to usual tableaux, $|\alpha|$ does not equal 
$|\gamma|$. Instead, $|\alpha|$ is the number of $k$-bounded hooks 
in $\gamma$.  This
distinction becomes natural through a correspondence between $k+1$-cores 
and $k$-bounded diagrams.  This bijection between $\CC^{k+1}$ and $\CP^k$
was defined in \cite{[LMcore]} by the map
$$
\kbnd\left(\gg\right) = (\lambda_1,\ldots,\lambda_\ell)
$$
where $\lambda_i$ is the number of cells with a $k$-bounded hook
in row $i$ of $\gamma$.  Note that the number of $k$-bounded hooks 
in $\gamma$ is $|\lambda|$.  The inverse map relies on constructing a 
certain ``$k$-skew diagram" $\lambda/^k=\gg/\rho$ from 
$\lambda$, and setting $\core(\lambda)=\gamma$.
These special skew diagrams are defined:

\begin{definition} \label{def5}
For $\lambda\in\mathcal P^k$, the
``$k$-skew diagram of $\lambda$" is the diagram $\lambda/^k$
where

(i) row $i$ has length $\lambda_i$ for $i=1,\ldots,\ell(\lambda)$

(ii) no cell of $\lambda/^k$ has hook-length exceeding $k$

(iii) all squares below $\lambda/^k$ have hook-length exceeding $k$.

\noindent
\end{definition}

A convenient algorithm for constructing the diagram of $\lambda/^k$ is 
given by successively attaching a row of length $\lambda_i$ to the bottom of 
$(\lambda_1,\dots,\lambda_{i-1})/^k$ in the leftmost position so that no 
hook-lengths exceeding $k$ are created.

\begin{example} \label{exskew}
Given $\lambda =(4,3,2,2,1,1)$ and $k=4$,
\begin{equation*}
\lambda = {\tiny{\tableau*[scY]{  \cr  \cr  & \cr &
\cr & & \cr & & & }}}
\quad
\implies
\quad
\lambda/^4 = {\tiny{\tableau*[scY]{  
\cr  
\cr  & 
\cr \bl &  & 
\cr  \bl &\bl & & & 
\cr \bl &\bl& \bl & \bl &\bl & & & &
}}}
\qquad\implies
\;
\core(\lambda)= {\tiny{\tableau*[scY]{  
\cr  
\cr  & 
\cr  &  & 
\cr  & & & & 
\cr  && & & & & & &
}}}
\end{equation*}
\end{example}

\medskip

The analogy with usual tableaux is now more apparent, and we let
$\mathcal T^k_{\aa}(\mu)$ denote the set of all $k$-tableaux of 
shape $\core(\mu)$ and $k$-weight $\aa$.  When the $k$-weight 
is $(1^n)$, a $k$-tableau is called {\it ``standard"}.  
The ``$k$-Kostka numbers" $K_{\mu\alpha}^{(k)}=|T^k_\aa(\mu)|$
satisfy a triangularity property \cite{[LMcore]} similar  to that
of the Kostka numbers: for $k$-bounded partitions $\lambda$ and $\mu$,
\begin{equation}
\label{trikostka}
K_{\mu\lambda}^{(k)}=0 \quad\text{when}\quad \mu \ntrianglerighteq
\lambda \quad
\text{ and }\quad K_{\mu\mu}^{(k)}=1\,.
\end{equation}

Given this triangularity, the inverse of 
$||K^{(k)}_{\mu\lambda}||_{\lambda,\mu\in \mathcal P^k}$ exists.  
Our main object of study can now be defined by
$||K^{(k)}||^{-1}$, denoted $||\bar K^{(k)}||$.

\begin{definition}
\label{kschurdef}
For any $\lambda\in\mathcal P^k$,
the ``$k$-Schur function" is defined
\begin{equation}
s_{\lambda}^{(k)} = \sum_{\mu \trianglerighteq \lambda} 
\bar K_{\mu \lambda}^{(k)} h_{\mu}\, .
\end{equation}
\end{definition}

A number of properties held by $k$-Schur functions suggest that these 
elements play the role of the Schur functions in the subspace $\Lambda^{k}$.  
First, the definition implies that the set 
$\left\{s_\lambda^{(k)}\right\}_{\lambda_1\leq k}$ 
forms a basis of $\Lambda^{k}$, and
that for any $\lambda\in\mathcal P^k$,
\begin{equation}
h_{\lambda}= \sum_{\mu \trianglerighteq \lambda} 
K_{\mu \lambda}^{(k)} s_{\mu}^{(k)}\, .
\end{equation}
In \cite{[LMkschur]} it was shown that these functions satisfy  
the  ``$k$-Pieri formula": for $\nu_1,\ell\leq k$, 
\begin{equation}
\label{kpieri}
h_\ell\,s_\nu^{(k)} = \sum_{\lambda\in H_{\nu,\ell}^{(k)}}
s_\lambda^{(k)}
\end{equation}
where the sum is over partitions of the form:
$$
H_{\nu,\ell}^{(k)}=
\left\{ \lambda \, \Big| \, 
\lambda/\nu={\rm horizontal~} \ell\text{-}{\rm strip} \quad {\rm and} \quad
\lambda^{\omega_k}/\nu^{\omega_k}={\rm vertical~} \ell\text{-}{\rm strip} 
\right\} \,.
$$
It was also shown more generally that if $K_{\nu/\mu,\lambda}^{(k)}$ is 
the number of skew tableaux of shape $\core(\nu)/\core(\mu)$ and 
$k$-weight $\lambda$, then
\begin{equation} \label{eqhons}
h_{\lambda} s_{\mu}^{(k)} = \sum_{\nu} K_{\nu/\mu,\lambda}^{(k)} \, 
s_{\nu}^{(k)} \, .
\end{equation}

\medskip

Another example of a Schur property held by $k$-Schur functions 
is drawn from the $\omega$-involution, defined as the homomorphism 
$\omega(h_i)=e_i$. In particular, $\omega$ maps a Schur function
$s_\lambda$ to it's conjugate $s_{\lambda'}$.   Using a refinement
of partition conjugation that arose in \cite{[LLM],[LMfil]}, 
it was shown in \cite{[LMkschur]} that
\begin{equation}
\label{omegainv}
\omega s_\lambda^{(k)} = s_{\lambda^{\omega_k}}^{(k)} \, ,
\end{equation}
where $\la^{\om_k} =  \kbnd\big(\core(\lambda)'\big)$
is the ``$k$-conjugate" of $\la$.  This result led to the
property:
\begin{equation}
\label{proposmall}
s_{\lambda}^{(k)}=s_{\lambda}\;\text{when}\;
h(\lambda) \leq k \,.
\end{equation}

\medskip

In the spirit of Schur function theory, it is conjectured that the
``$k$-Littlewood-Richardson coefficients" in
\begin{equation}
\label{klr}
s_\lambda^{(k)} s_\mu^{(k)} =\sum_{\nu: \nu_1\leq k}
c_{\lambda\mu}^{\nu, k} s_\nu^{(k)}\,,
\end{equation}
are positive numbers.  Our development here will prove that
in certain cases, these coefficients are the Gromov-Witten
invariants thus proving positivity in these cases.  Note that
given the action of the $\omega$ involution on $k$-Schur functions,
the $k$-Littlewood-Richardson coefficients satisfy
\begin{equation}\label{symmetry}
c_{\lambda \mu}^{\nu, k} = c_{\lambda^{\omega_k} 
\mu^{\omega_k}}^{\nu^{\omega_k}, k} \, . 
\end{equation}

\section{Hecke algebras, fusion rules, and the $k$-Schur functions}

Presented in \cite{[GW]} are generalized Littlewood-Richardson coefficients 
for $(\ell,n)$-representations of the Hecke algebras $H_{\infty}(q)$, when 
$q$ is an $n$-th root of unity.  These coefficients are
equivalent to the structure 
constants for the Verlinde (fusion) algebra associated to the 
$\widehat{su}(\ell)$-Wess-Zumino-Witten conformal field 
theories at level $n-\ell$.  In this section, we will use the $k$-Pieri rule 
to establish that for $k=n-1$, the $k$-Littlewood-Richardson 
coefficients contain these constants as special cases.

\medskip

\subsection{The connection}
From \cite{[GW]}, we recall a simple interpretation for
these ``$(\ell,n)$-Littlewood-Richardson coefficients" 
given in the language of symmetric functions.
For $n>\ell\geq 1$, consider the quotient 
$R^{\ell n} = \Lambda^{\ell}/I^{\ell n}$
where $I^{\ell n}$ is the ideal generated by
Schur functions that have exactly $n-\ell+1$
rows of length smaller than $\ell$:
$$
I^{\ell n}= \Bigl \langle s_{\lambda} \, 
\Big| \, 
 \# \{ j \, |  \, \lambda_j < \ell\}=n-\ell+1 \Bigr \rangle \, .
$$
A basis for $R^{\ell n}$ is given by the set 
$\{ s_{\lambda}\}_{\lambda \in \Pi^{\ell n}}$ where the indices are
partitions in:
$$
\Pi^{\ell n} = \{\lambda\in \mathcal P :
\lambda_1\leq \ell \;\text{and}\;
 \# \{ j \, | \, \lambda_j < \ell\}\leq n-\ell 
\}
\,.
$$
The $(\ell,n)$-Littlewood-Richardson coefficients of interest here are
simply $a^\nu_{\lambda \mu}$ in
\begin{equation}
\label{gwpieri}
s_\lambda s_\mu = \sum_\nu a^\nu_{\lambda \mu} \, s_\nu \mod  I^{\ell n}
 \, ,\quad
\text{where $\lambda,\mu,\nu \in  \Pi^{\ell n}$. }
\end{equation}
It is in this context that we prove the coefficients $a_{\lambda\mu}^{\nu}$
are none other than
$k$-Littlewood-Richardson coefficients when $k=n-1$.

\begin{remark}  
The results of \cite{[GW]} are presented in a transposed form,
where they instead work with the ideal
$\langle s_{\lambda} \, | \, \lambda_1-\lambda_{\ell}=n-\ell+1 \rangle$ in
$\mathbb Z[e_1,\dots,e_l]$. Their $(\ell,n)$-Littlewood-Richardson 
coefficients $d_{\lambda \mu}^{\nu}$ are our $a_{\lambda' \mu'}^{\nu'}$, 
for $\lambda',\mu',\nu' \in  \Pi^{\ell n}$.
\end{remark}

To provide some insight into how this connection arose, consider
the special case of Eq.~\eqref{gwpieri} with $\lambda=(1)$:
\begin{equation} \label{eqschur}
s_1 s_{\mu}=
\sum_{\nu \, : \, \mu\subset\nu\in \Pi^{\ell n}\atop 
|\nu|=|\mu|+1} s_{\nu} \mod  I^{\ell n}  \, ,
\end{equation}
and define a poset by letting $\mu\prec\!\!\cdot\nu$ for all $\nu$ in 
the summand.  Frank Sottile brought this poset to our attention and asked if 
it was related to our study \cite{[LMcore]} of the $k$-Young lattice $Y^k$.  
$Y^k$  is defined by the $k$-Pieri rule, where $\mu<\!\!\cdot\nu$ when 
$\nu\in H_{\mu,1}^{(k)}$.  Investigating his question, we discovered the posets 
can be connected through the principal order ideal $L^{k}(\ell,m)$ generated by 
an  $\ell\times m$ rectangle in $Y^k$.  In \cite{[LMrec]}, we found 
that the vertices of $L^{k}(\ell,m)$ are the partitions contained in an 
$\ell \times m$ rectangle with no more than $k-\ell+1$ rows shorter than $k$, and 
that $\mu$ covers $\lambda$ in this poset if and only if $\lambda \subseteq \mu$ 
and $|\lambda|+1=|\mu|$.  Therefore, the elements of $L^{k}(\ell,\infty)$ are 
precisely those of $\Pi^{\ell n}$ (given $k=n-1$).  Since the $k$-Young lattice 
was defined by multiplication by $s_1$, we have
\begin{equation} 
\label{eqkschur}
s_1 s^{(k)}_{\lambda} = 
\sum_{\mu \, : \, \lambda\subset \mu \in  \Pi^{\ell n}\atop
|\mu|=|\lambda|+1}  s^{(k)}_{\mu} + {\rm ~other~terms} \, ,
\end{equation}
where ``other terms'' are $k$-Schur functions indexed by 
$\mu\not\in\Pi^{\ell n}$.  The likeness of \eqref{eqschur} 
and \eqref{eqkschur} led us to surmise the following result:

\begin{theorem}
\label{modth}
For any partition $\lambda\in\mathcal P^{n-1}$,
\begin{equation}
s_{\lambda}^{(n-1)} \mod {I^{\ell n}}= 
\begin{cases}
s_{\lambda} & \text{if } \lambda \in \Pi^{\ell n} \\
0 & {\rm otherwise}
\end{cases}
\end{equation}
\end{theorem}

Before proving this theorem, we mention several implications.
Since all partitions in $\Pi^{\ell n}$ are $(n-1)$-bounded 
($\lambda_1\leq \ell \leq n-1$), the set of $k$-Schur functions 
indexed by partitions in $\Pi^{\ell n}$ forms a natural basis for the
quotient $R^{\ell n}$.  Computation modulo the ideal 
$I^{\ell n}$ is trivial in this basis.  In particular,  
the structure constants under consideration are simply
certain $k$-Littlewood Richardson coefficients in \eqref{klr}.

\begin{corollary}
For all $\lambda,\mu,\nu \in  \Pi^{\ell n}$,
$$a_{\lambda \mu}^{\, \nu}= c_{\lambda \mu}^{\nu, n-1}\,.$$
\end{corollary}

Another consequence of our theorem produces a tableau interpretation for
the dimension of the representations $\pi_{\lambda}^{(\ell,n)}$, for
$\lambda' \in \Pi^{\ell n}$, of the 
Hecke algebras $H_{\infty}(q)$, when $q$ is an $n$-th root of unity
(see \cite{[GW]} for details on these representations).

\begin{corollary}
For $\lambda' \in \Pi^{\ell n}$, the dimension of the representation 
$\pi^{(\ell,n)}_{\lambda}$ is the number of standard $(n-1)$-tableaux of shape 
$\core(\lambda')$\footnote{Equivalently, this is the number of reduced words for 
a certain affine permutation $\sigma_{\lambda'}\in\hat S_{n}/S_{n}$.  
See \cite{[LMcore]} for the precise correspondence.}.
\end{corollary}
\begin{proof}
Let $m=|\lambda|$, and $k=n-1$.  In \cite{[GW]}, it is shown that
the dimension of $\pi^{(\ell,n)}_{\lambda}$ is the coefficient of 
$s_{\lambda'}$ in $s_1^m\mod I^{\ell n}$.  By Theorem~\ref{modth},
this is the coefficient of $s_{\lambda'}^{(k)}$ in 
the $k$-Schur expansion of $s_1^{m}=h_{1^m}$.  Using 
Definition~\ref{kschurdef} for $k$-Schur functions,
this coefficient is $K_{\lambda' 1^m}^{(k)}$, or the number of 
standard $k$-tableaux of shape $\core(\lambda')$.  
\end{proof}

The Verlinde (fusion) algebra  of the Wess-Zumino-Witten model associated 
to $\widehat{su}(\ell)$ at level $n-\ell$ is isomorphic to the quotient of
$R^{\ell n}$ modulo the single relation $s_{\ell}\equiv 1$ 
\cite{[Kac],[Wa],[GW]}.  The fusion coefficient 
$\mathcal N_{\lambda \mu}^{\nu}$ is defined for
$\lambda',\mu',\nu'\in\mathcal P^{\ell-1\,,n-1}$ by
$$
L(\lambda)\otimes_{n-\ell} L(\mu) = 
\oplus
\mathcal N_{\lambda \mu}^{\nu}
L(\nu)
\,,
$$
where the fusion product $\otimes_{n-\ell}$ is the reduction of the tensor
product of integrable representations with highest weight $\lambda$ and $\mu$
via the representation at level $n-\ell$ of $\widehat{su}(\ell)$.
Thus, our results imply that

\begin{corollary}
For all $\lambda,\mu,\nu$ inside an $(n-\ell)\times(\ell-1)$ rectangle,
$$\mathcal N_{\lambda \mu}^{\nu}= c_{\lambda' \mu'}^{\hat\nu, n-1}\,,
$$ 
where $\hat\nu=(\ell^{(|\lambda|+|\mu|-|\nu|)/\ell},\nu')$.
\end{corollary}

\subsection{Proof of the connection}

To prove Theorem~\ref{modth}, we use two preliminary properties.
For simplicity, since 
$\Lambda^{\ell}= \Lambda/  \langle h_{\ell+1},h_{\ell+2},\dots  \rangle$,
we will instead work with the ideal $\mathcal I$ in $\Lambda$, where 
$$
\mathcal I = \bigl \langle
s_{\lambda} \, \big | \, 
 \# \{ j \, |  \, \lambda_j < \ell\}=n-\ell+1  \bigr \rangle \cup 
\bigl \langle  h_i \,\big |\, i>\ell \bigr \rangle
\, ,
$$
and in the remainder of this section, $k$ will always stand for $n-1$.

\begin{property}
\label{biggerm}
For any $k$-bounded partition $\lambda$ and $\ell\leq k$,
$s_\lambda^{(k)}\equiv_{\mathcal I} 0$ when $\lambda_1>\ell$.
\end{property}
\begin{proof}
Since $\mu\geq\lambda$ implies that $\mu_1\geq\lambda_1$,
the unitriangular relation between $\{s_\lambda^{(k)}\}$
and $\{h_\lambda\}$ implies
$$
s_{\lambda}^{(k)} 
= \sum_{\mu :\mu_1>\ell} * \, h_{\mu} 
\, .
$$
The claim thus follows since $h_\mu \in \mathcal I$ when $\mu_1>\ell$.
\end{proof}

\begin{property}
\label{addms}
For any $k$-bounded partition $\lambda$ with $\lambda_1\leq \ell$,
\begin{equation}
s_\lambda^{(k)}\equiv_{\mathcal I} 0 \implies
s_{(\ell^m,\lambda)}^{(k)} \equiv_{\mathcal I} 0\quad\text{for all} \;\;
m\geq 0\,.
\end{equation}
\end{property}

\begin{proof}
The $k$-Pieri rule \eqref{kpieri} implies in particular, that any $k$-Schur occurring in 
the expansion of $h_\ell s_\nu^{(k)}$ is indexed by a partition obtained by 
adding a horizontal $\ell$-strip to $\nu$.  Thus, when $\ell\geq \nu_1$, we have
\begin{equation}
\label{unipieri}
h_\ell s_\nu^{(k)} =
s_{(\ell,\nu)}^{(k)}\,+\,
\sum_{\mu:\mu_1>\ell
\atop
\mu\in H_{\nu,\ell}^k} s_\mu^{(k)}\,.
\end{equation}
Starting from $s_\lambda^{(k)}\equiv_{\mathcal I} 0$, and assuming by 
induction that $s_{(\ell^{m-1},\lambda)}^{(k)}\equiv_{\mathcal I} 0$,
the claim follows from Property~\ref{biggerm} and the previous
expression \eqref{unipieri},
$$
0\;\equiv_{\mathcal I}\;
h_\ell\,s_{(\ell^{m-1},\lambda)}^{(k)}
\;= \;
s_{(\ell^m,\lambda)}^{(k)}+ 
 \sum_{\gamma  \, : \, \gamma_1>\ell}
* \, s_\gamma^{(k)}
\;\equiv_{\mathcal I}\;
s_{(\ell^m,\lambda)}^{(k)} \,.
$$
\end{proof}

\subsection{Proof of Theorem~\ref{modth}}
Recall $n=k+1$, and that
$\lambda\in \Pi^{\ell,k+1}$ has the form $\lambda=(\ell^m,\mu)$ for some 
$\mu\in\mathcal P^{\ell-1\, k}$.  First, by induction on $m$ we prove
that $s_\lambda^{(k)}\equiv_{\mathcal I} s_\lambda$ for each such $\lambda$.
Since $h(\lambda)\leq k$ when $m=0$, $s_\lambda^{(k)}=s_\lambda$ by 
\eqref{proposmall}.  By induction, assuming
$s_{(\ell^m,\mu)}^{(k)}\equiv_{\mathcal I} s_{(\ell^m,\mu)}$, we have
$h_\ell\,s_{(\ell^m,\mu)}^{(k)} \equiv_{\mathcal I} h_\ell\,s_{(\ell^m,\mu)}$.
On the other hand, since $s_\gamma\equiv_{\mathcal I} 0$ when $\gamma_1>\ell$, 
Identity \eqref{unipieri} implies 
$h_\ell\,s_{(\ell^m,\mu)}^{(k)}
\equiv_{\mathcal I}s_{(\ell^{m+1},\mu)}^{(k)}$.
Therefore,
$$
h_\ell\,s_{(\ell^m,\mu)}
\equiv_{\mathcal I}s_{(\ell^{m+1},\mu)}^{(k)}\,.
$$
The claim then follows by noting that the Pieri rule gives an
expansion similar to \eqref{unipieri} for $h_\ell\,s_{\ell^m,\mu}$, 
implying that $h_\ell\,s_{(\ell^m,\mu)} \equiv_{\mathcal I}s_{(\ell^{m+1},\mu)}$.

It remains to prove that $s_\eta^{(k)}\equiv_{\mathcal I} 0$ when 
$\eta\not\in\Pi^{\ell,k+1}$.  Since Property~\ref{biggerm} proves 
the case when $\eta_1>\ell$, we must show
$s_{\eta}^{(k)}\equiv_{\mathcal I} 0$ for any $\eta$ in the set:
$$ \mathcal Q=\bigl\{(\ell^m,\beta)\in\mathcal P :
\beta_1<\ell\;\text{and}\; \ell(\beta)\geq k-\ell+2\bigr\}\,.
$$  
Our proof is inductive, using an order defined on $\mathcal Q$
as follows: $\eta=(\ell^a,\beta)\preceq (\ell^b,\alpha)=\mu$ if
$\ell(\beta)<\ell(\alpha)$ or if $\ell(\beta)=\ell(\alpha)$ and 
$\eta\unrhd\mu$ (this is a well-ordering if we restrict ourselves 
to $|\mu|=|\eta|$).  Our base case includes partitions $\eta=(\ell^a,\beta)$ 
with $\beta_1<\ell$ and $\ell(\beta)=k-\ell+2$.  In this case, $h(\beta)\leq k$ 
implies $s_\beta^{(k)}=s_\beta$ from \eqref{proposmall}, and since
$s_\beta\in {\mathcal I}$ when $\beta$ has $k-\ell+2$ parts 
smaller than $\ell$, we have $s_\beta^{(k)}\equiv_{\mathcal I} 0$.
Property~\ref{addms} then proves
$s_\eta^{(k)}\equiv_{\mathcal I} 0$ in this case.

Now assume by induction that $s_{\eta}^{(k)} \equiv_{\mathcal I} 0$ for 
all $\eta\in \mathcal Q$ such that $\eta\prec\mu$, where
$\mu=(\ell^b,\alpha)$ with $\alpha_1 < \ell$ and 
$\ell(\alpha)>k-\ell+2$.  With $r<\ell$ denoting the last part 
of $\mu$ (and thus also the last part of $\alpha$), 
let $\mu=(\hat\mu,r)=(\ell^b,\hat\alpha,r)$
and note that $\hat\mu\prec\mu$.    Thus, using
the induction hypothesis and the $k$-Pieri rule, we have 
$$
0\,\equiv_{\mathcal I}\,
s_r s_{\hat\mu}^{(k)} = s_{\mu}^{(k)}+ 
\sum_{\nu \in
H_{\hat\mu,r}^{(k)}\setminus \{ \mu\}
} s_{\nu}^{(k)} \, ,
$$
and it suffices 
to show that $s_\nu^{(k)}\equiv_{\mathcal I} 0$ for all $\nu\in 
H_{\hat\mu,r}^{(k)}\setminus \{ \mu\}$.  Property~\ref{biggerm} 
proves this immediately for any $\nu$ with $\nu_1>\ell$, and thus
we shall consider only $\ell$-bounded $\nu$.  Two properties of such 
$\nu$ follow since $\nu$ is obtained by adding a horizontal $r$-strip to 
$\hat\mu=(\ell^b,\hat\alpha)$:  $\nu\rhd\mu$, and
$\nu=(\ell^b,\beta)$, where $\ell(\beta)\leq\ell(\hat\alpha)+1=\ell(\alpha)$.  
Thus, if these $\nu$ lie in $\mathcal Q$, then $\nu \prec \mu$
and our claim follows from the induction hypothesis.
Since each such $\nu$ is obtained by adding a horizontal strip 
to $\hat\mu=(\ell^b,\hat\alpha)$, and $\ell(\alpha)>k-\ell+2$, we have
$\ell(\beta)\geq\ell(\hat\alpha)\geq k-\ell+2$.  Thus, these 
$\nu=(\ell^b,\beta)$ all lie in $\mathcal Q$ except in the 
case that $\ell(\beta)=\ell(\hat\alpha)=k-\ell+2$ and $\beta_1=\ell$.
The following paragraph explains why, in this case, $\nu\not\in
H_{\hat \mu,r}^{(k)}$, and thus never arise.

Given $\ell(\beta)=k-\ell+2$ and $\beta_1=\ell$, 
$h(\beta)>k$ implies there is no cell in position 
$X=(1,\ell(\nu)-\ell(\beta))$ of $\nu/^k$.  Assume by contradiction that
$\nu\in H_{\hat\mu,r}^{(k)}$ -- hence, in particular, that
$\nu^{\omega_k}/\hat\mu^{\omega_k}$ is a vertical strip.  
Since $\ell(\hat\alpha)=\ell(\beta)=k-\ell+2$ and $\hat\alpha_1<\ell$ imply
$h(\hat\alpha)\leq k$, there is a cell in $\hat\mu/^k$ in position
$(1,\ell(\mu)-\ell(\hat \alpha))=X$.  Since the height of $\nu/^k$ and
$\hat \mu/^k$ are equal, but position $X$ 
is empty in $\nu/^k$ and filled in $\hat \mu/^k$,
the first column of $\nu/^k$ is shorter than that of
$\hat\mu/^k$, implying $\nu^{\omega_k}/{\hat \mu}^{\omega_k}$
is not a vertical strip.  By contradiction, $\nu \not \in  H_{\hat\mu,r}^{(k)}$
as claimed.  Here is an example with
$\hat \mu= (4,4,2,2,1,1)$, $\nu= (4,4,4,2,1,1)$, $n=4$ and $k=6$:
$$
\hat \mu/^k=
{\tiny{\tableau[scY]{| |,|X ,|\bl,,,,|\bl,\bl,,,,
}}} 
\quad \quad \quad \nu/^k = 
{\tiny{\tableau[scY]{| |,|\bl X ,,,,|\bl,,,,|\bl,\bl,,,,
}}}
$$
\hfill$\square$

\medskip

\section{Quantum cohomology}

Witten \cite{[Wi]} proved that the Verlinde algebra of $\widehat{u}(\ell)$
at level $n-\ell$ and the quantum cohomology of the Grassmannian $Gr_{\ell n}$ 
are isomorphic (see also \cite{[Ag]}).  Since $u(\ell)=su(\ell) \times u(1)$,
the connection between $k$-Schur functions and the fusion coefficients of 
$\widehat{su}(\ell)$ at level $n-\ell$ given in the last section implies 
that there is also a connection between $k$-Schur functions and the quantum 
cohomology of the Grassmannian.  We now set out to make this connection explicit.

\medskip

Recall from the introduction that the quantum structure constants,
or 3-point Gromov-Witten invariants $C_{\lambda\mu}^{\nu,d}$,
arise in the expansion, for $\lambda,\mu\in  \mathcal P^{\ell n}$, 
\begin{equation} 
\label{prodss}
s_\lambda \, s_{\mu} 
= \sum_{d\geq 0,\,
\nu\in\mathcal P^{\ell n}\atop
|\nu|=|\lambda|+|\mu|-d n}
q^d\, C_{\lambda\mu}^{\nu,d}\, s_\nu 
\mod J^{\ell n}_q 
\, ,
\end{equation}
where
$$
J^{\ell n}_q=\langle e_{n-\ell+1},\ldots,e_{n-1},e_n+(-1)^{\ell}q\rangle
\,.
$$
Our main goal is to prove that the $k$-Schur function basis gives
a direct route to these constants. In particular, by 
determining the value of a $k$-Schur function modulo this ideal,
we will see that the Gromov-Witten invariants arise as 
special cases of the $k$-Littlewood-Richardson coefficients.

\medskip
Again for simplicity, we work in $\Lambda/\mathcal J_q$, 
where $\mathcal J_q$ is the ideal 
$$
\mathcal J_q = J^{\ell n}_q \cup \langle
 h_i\, | \,  i>\ell
\rangle
\,.
$$
Theorem~\ref{modth} reveals that a $k$-Schur function modulo the ideal 
$\mathcal I$ is a Schur function when $\lambda\in\Pi^{\ell n}$ and is 
otherwise zero.  By showing that $\mathcal I$ is a subideal of $\mathcal J_q$,
our task to determine a $k$-Schur function mod $\mathcal J_q$ is thus
reduced to examining what happens to a usual Schur function $s_\lambda$ 
mod $\mathcal J_q$ in the special case that $\lambda\in\Pi^{\ell n}$.

\begin{proposition} \label{subideal}
If $f \in \mathcal I$, then $f \in \mathcal J_q$.
\end{proposition}
\begin{proof}
It suffice to prove that $s_\lambda \in \mathcal J_q$
when $\lambda=(\ell^m,\alpha)$, for some $m$ and partition $\alpha$
such that $\alpha_1<\ell $ and $\ell(\alpha)=n-\ell+1$.  When $m=0$, 
the result follows from the Jacobi-Trudi determinantal formula since 
the first row of the determinant of $s_{\alpha}$ has entries
$e_{n-\ell+1},\dots,e_{n+\alpha_1-\ell}\in \mathcal J_q$ 
given $\alpha_1<\ell$.  Assuming by induction that
$s_{(\ell^{m},\alpha)}\in \mathcal J_q$, since the Pieri rule implies 
$h_{\ell} s_{(\ell^m,\alpha)} = s_{(\ell^{m+1},\alpha)} 
\mod \langle h_{\ell+1},h_{\ell+2},\dots \rangle$, the
result follows by induction.
\end{proof}

Now to determine the value of a usual Schur function $s_\lambda$ mod 
$\mathcal J_q$ for partitions in $\Pi^{\ell n}$, we shall use
an important result from \cite{[BCF]}, where the theory of rim-hooks 
was used to study the Schur functions modulo $\mathcal J_q$.
To state their result, we first recall the necessary definitions.  
An {\it``$n$-rim hook"} is a connected skew diagram of size $n$
that contains no $2 \times 2$ rectangle.  ``$\rcore(\lambda)$"
denotes the $n$-core of $\lambda$, obtained by removing as many 
$n$-rim hooks as possible from the diagram of $\lambda$ (this
is well-defined since the order in which rims are removed is known 
to be irrelevant \cite{[JK]}).  The width of a rim hook is the number 
of columns it occupies minus one.  Given a partition $\lambda$, 
let ``$d_{\lambda}$" be the number of $n$-rim hooks that are removed 
to obtain $\rcore(\lambda)$.  Also, let ``$\epsilon_{\lambda}$" 
equal $d_{\lambda}(\ell-1)$ minus the sum of the widths of these 
rim hooks.  This given, in \cite{[BCF]} Eq.~(19) it is shown for
$\lambda \in \mathcal P^{\ell}$, that
\begin{equation} \label{smodJ}
s_{\lambda} \mod J^{\ell n}_q= \left\{
\begin{array}{ll}
(-1)^{\epsilon_{\lambda}} 
q^{d_{\lambda}} s_{\rcore(\lambda)} & {\rm if~} \rcore(\lambda)
\in \mathcal P^{\ell n} \\
0 & {\rm otherwise}
\end{array}\right.
\end{equation}

\medskip

This result helps us prove that in the special case that
$\lambda\in\Pi^{\ell n}$, $s_{\lambda}\equiv_{\mathcal J_q}
q^{d_{\lambda}} s_\nu$  for a partition $\nu$ 
obtained using the following operators:

\begin{definition}
For $\lambda\in\Pi^{\ell n}$, $\blacktriangle(\lambda)$ is the partition
obtained by adding an $n$-rim hook to $\lambda$ starting in column $\ell$ 
and ending in the first column.  For $\lambda\in\Pi^{\ell n}$ that is not 
an $n$-core, $\blacktriangledown(\lambda)$ is the partition obtained by 
removing an $n$-rim hook from $\lambda$ starting in the first column 
of $\lambda$.
\end{definition}

Note that $\blacktriangledown$ is well-defined since
when $\lambda \in \Pi^{\ell n}$ is not an $n$-core, $\ell(\lambda)>n-\ell$.
Thus for $r=\ell(\lambda)-(n-\ell)$, $\lambda_r=\ell$ and 
$h_{(r,1)}(\lambda)=n$ implying an $n$-rim hook can be removed
starting in the first column of $\lambda$ and ending in the last column
$\ell$.  Since the difference between the heights of the starting point 
and the ending point is $n-\ell$,
\begin{equation}
\label{minusremain}
\blacktriangledown : \left\{\lambda \, | \, \lambda\in\Pi^{\ell n}\,\, \&\,\, \lambda\neq
\mbox{$n$-core}\right\} \; \to\;\Pi^{\ell n}
\,.
\end{equation}
Similarly, for any $\lambda \in \Pi^{\ell n}$, the difference 
between the heights of the first column and column $\ell$ is at 
most $n-\ell$.  Thus, an $n$-rim hook can be added to $\lambda$ starting 
from column $\ell$ and ending in the first column.  Since the difference 
in heights of the starting point and ending point of the added $n$-rim hook 
is $n-\ell$, we have that
\begin{equation}
\label{plusremain}
\blacktriangle : \Pi^{\ell n}\to\Pi^{\ell n}
\,.
\end{equation}
By construction,
as long as $\blacktriangledown(\lambda)$ is defined, we have
\begin{equation}
\label{plusminus}
\blacktriangledown(\blacktriangle(\lambda))=\lambda\; \quad {\rm and}\; \quad
\blacktriangle(\blacktriangledown(\lambda))=\lambda
\,.
\end{equation}

\begin{proposition} \label{propmod}
For $\lambda \in \Pi^{\ell n}$,
\begin{equation} \label{eqpropo}
s_{\lambda}\equiv 
q^{d_{\lambda}} \, s_\nu
\mod {\mathcal J}_q 
\, ,
\end{equation}
where $\nu={\blacktriangledown^{d_\lambda}(\lambda)}\in\mathcal P^{\ell n}$.
\end{proposition}
\begin{proof} 
We have $s_{\lambda} \equiv_{{\mathcal J}_q }
(-1)^{\epsilon_\lambda}q^{d_{\lambda}} \, s_{\rcore(\lambda)}$
by \eqref{smodJ}.
When $\lambda$ is an $n$-core, then $\lambda \in \Pi^{\ell n}$ implies
$\lambda \in \mathcal P_{\ell n}$.  Thus $\rcore(\lambda)=\lambda$
and \eqref{eqpropo} holds with $d_\lambda=0$.  Otherwise, $\rcore(\lambda)$ 
is obtained by removing $d_\lambda$ $n$-rim hooks in any order.
Thus, by successively applying $\blacktriangledown$, we obtain
$\rcore(\lambda)=\blacktriangledown^{d_\lambda}(\lambda)$. 
Since $\blacktriangledown$ preserves $\Pi^{\ell n}$ by \eqref{minusremain},
$\rcore(\lambda)\in\mathcal P^{\ell n}$.  Further,
$\epsilon_{\lambda}=d_{\lambda}(\ell-1)-d_{\lambda}(\ell-1)=0$ since 
each removed $n$-rim hook has width $\ell-1$. 
\end{proof}

In this notation, we can now determine the value of a 
$k$-Schur function mod $\mathcal J_q$.

\begin{theorem} \label{quantummod}
For any $k$-bounded partition $\lambda$,
$$
s_\lambda^{(n-1)} \mod {\mathcal J_q} =
\begin{cases}
q^{d_{\lambda}} s_{\nu}  
& {\rm if~} \lambda \in \Pi^{\ell n} \\
0 & {\rm otherwise}
\end{cases}
$$
where 
$\nu=\rcore(\lambda)={\blacktriangledown^{d_\lambda}(\lambda)}\in\mathcal 
P^{\ell n}$.
\end{theorem}
\begin{proof}  
Proposition~\ref{subideal} gives that
$\mathcal I$ is a subideal of $\mathcal J_q$, implying
$$
s_{\lambda}^{(n-1)} \mod {\mathcal J_q} = 
\left(s_{\lambda}^{(n-1)} \mod {\mathcal I} \right)
\mod {\mathcal J_q} \, .
$$
For $\lambda\not\in\Pi^{\ell n}$,
$s_{\lambda}^{(n-1)} \mod {\mathcal I} = 0$ by 
Theorem~\ref{modth}.  For $\lambda\in\Pi^{\ell n}$, Theorem~\ref{modth} 
implies that $s_{\lambda}^{(n-1)} \mod {\mathcal I} = s_{\lambda}$, and
the claim then follows by further moding out by $\mathcal J_q$
according to Proposition~\ref{propmod}.
\end{proof}

This theorem enables us to connect the quantum product to the product
of $k$-Schur functions.

\begin{theorem}  \label{gromovwitten}
For $\lambda,\mu,\nu \in \mathcal P^{\ell n}$,
the 3-point Gromov-Witten invariants
$C_{\lambda \mu}^{\nu,d}$ are 
\begin{equation}
C_{\lambda \mu}^{\nu,d} = c_{\lambda \mu}^{\hat \nu ,n-1} \, ,
\end{equation}
where $\hat \nu = \blacktriangle^d(\nu)$, and where 
$c_{\lambda \mu}^{\hat \nu ,n-1}$ is a $k$-Littlewood-Richardson
coefficient.
\end{theorem}
\begin{proof} 
For $\lambda,\mu\in\mathcal P^{\ell n}$, 
\eqref{prodss} shows that $C_{\lambda \mu}^{\nu,d}$ arise in
the expansion
\begin{equation}
\label{side1}
s_{\lambda} \, s_{\mu}
  \equiv  
\sum_{d\geq 0,\,
\nu\in\mathcal P^{\ell n}\atop
|\nu|=|\lambda|+|\mu|-d n} C_{\lambda \mu}^{\nu,d} \, q^d \, s_\nu 
 \mod {\mathcal J_q} 
\,.
\end{equation}
On the other hand, since $\lambda,\mu\in\mathcal P^{\ell n}$ have
hook-length smaller than $n$, \eqref{proposmall} implies
that $ s_\lambda^{(n-1)} s_\mu^{(n-1)} = s_\lambda s_\mu $.
Therefore, applying Proposition~\ref{propmod} to the
$k$-Schur expansion of this product gives
\begin{equation}
\label{side2}
s_{\lambda}\, s_{\mu}
= \sum_{\gamma:
|\gamma|=|\lambda|+|\mu|} 
c_{\lambda \mu}^{\gamma, n-1} s_{\gamma}^{(n-1)} 
\equiv \sum_{\gamma \in \Pi^{\ell n}} c_{\lambda \mu}^{\gamma, n-1} 
\, q^{d_{\gamma}}\, s_{\beta}
\mod {\mathcal J_q} 
\,,
\end{equation}
where $\beta= \blacktriangledown^{d_\gamma}(\gamma)$.
Taking the coefficient of $q^d  s_\nu$ in \eqref{side1} and
\eqref{side2} implies
$$
C_{\lambda\mu}^{\nu,d}  =
\sum_{\gamma \in \Pi^{\ell n}
\atop \gamma \, : \, \nu=\blacktriangledown^{d}(\gamma) } 
c_{\lambda \mu}^{\gamma, n-1} \, .
$$
Since $\nu=\blacktriangledown^{d}(\gamma)\in\mathcal P^{\ell n}\subseteq
\Pi^{\ell n}$, we can apply $\blacktriangle$ to find there is a unique
$\gamma$ in the right summand.  That is, $\blacktriangle^{d}(\nu)=\gamma$
by \eqref{plusminus}.
\end{proof}

\medskip

It is important to note that the quantum structure constants
$C_{\lambda\mu}^{\nu,d}$ are indexed by $\lambda,\mu,\nu\in\mathcal P^{\ell n}$.
We have now seen that these numbers are precisely $k$-Littlewood-Richardson 
coefficients in the relevant cases.  However, since there are far more 
$k$-Littlewood-Richardson coefficients than Gromov-Witten invariants we 
naturally sought the larger geometric picture that would be explained by 
the complete set of $k$-Littlewood Richardson coefficients.  In discussions 
with Mark Shimozono about this problem, he conjectured that the $k$-Schur 
functions form the Schubert basis for the homology of affine (loop) Grassmannian
of $GL_{k+1}$, and that the expansion coefficients of the coproduct of $k$-Schur functions in 
terms of $k$-Schur functions gives the integral cohomology of loop Grassmannian
(see e.g. \cite{[Bot],[KK]} for more details on the loop Grassmannian).  This 
conjecture has been checked extensively with the
computer.  Corollary~\ref{gromovwitten} provides further evidence for this assertion 
based on the existence \cite{[DP]} of a surjective ring homomorphism from 
the homology of the loop Grassmannian onto the quantum cohomology of the 
Grassmannian at $q=1$.

\section{Dual $k$-schur functions}

While the homology of the loop Grassmannian is isomorphic to $\Lambda^k$,
the cohomology is isomorphic to $\Lambda/\mathfrak J^{(k)}$ for the ideal
$$
\mathfrak J^{(k)}=\langle m_\lambda:\lambda_1 >  k\rangle\,.
$$  
The interplay between homology and cohomology suggests that there
is a fundamental basis for $\Lambda/\mathfrak J^{(k)}$ that is
closely tied to the $k$-Schur basis.  Here, we introduce a family of
functions defined by the $k$-weight of $k$-tableaux and derive
a number of properties including a duality relation to the $k$-Schur 
functions.  In particular, it will develop that if the coproduct of
$k$-Schur functions in terms of $k$-Schur functions
indeed gives the integral cohomology of the loop Grassmannian, 
then these functions are the Schubert basis for the cohomology of the
loop Grassmannian.

\medskip

Recall that a Schur function can be defined as
\begin{equation}
s_{\lambda}= \sum_{T} x^{T} \, ,
\end{equation}
where the sum is over all tableaux of shape $\lambda$.  Several 
weight-permuting involutions have been defined on the set of tableaux
such as the Bender-Knuth involution \cite{[BK]}. 
From this, the Schur functions are combinatorially proven to be 
symmetric.  We extend these ideas by considering the family of functions 
that arises similarly from the set of $k$-tableaux (defined in 
\S~\ref{sect3} with $k$-weight).

\begin{definition} For any $\lambda\in\mathcal P^k$, 
the ``dual $k$-Schur function" is defined by 
\begin{equation}
\mathfrak S_{\lambda}^{(k)} = \sum_{T} x^{T} \, ,
\end{equation}
where the sum is over all $k$-tableaux of shape $\core(\lambda)$, and
$x^{T}=x^{k\text{-weight}(T)}$.
\end{definition}

An involution $\tau_a$ on the set of all $k$-tableaux generalizing the 
Bender-Knuth involution was given in \cite{[LMkschur]}.  This map
has the property that
\begin{equation}
\label{tauinvo}
\tau_a:\mathcal T_\alpha^k(\lambda) \to
\mathcal T_{\hat\alpha}^k(\lambda) \,,
\end{equation}
where $\hat \alpha=(\ldots,\alpha_{a+1},\alpha_a,\ldots,)$
is obtained by switching $\alpha_{a}$ and $\alpha_{b}$ in $\alpha$.
Since for any composition $\gamma$ (allowing zeros), the coefficient of 
$x^\gamma$ in  $\mathfrak S_{\lambda}^{(k)}$ is the number of $k$-tableaux 
of $k$-weight $\gamma$, this involution immediately implies the
symmetry of the dual $k$-Schur functions.

\medskip

\begin{proposition}  
\label{dualsym}
For any $k$-bounded partition $\lambda$,
$\mathfrak S_{\lambda}^{(k)}$ is a symmetric function.
\end{proposition}

\medskip

Since the $k$-Kostka number $K_{\lambda\alpha}^{(k)}$ denotes
the number of $k$-tableaux of shape $\core(\lambda)$
and $k$-weight $\alpha$, the symmetry of dual $k$-Schur functions
implies that
$$
\mathfrak S_\lambda^{(k)} = \sum_{\mu} K_{\lambda\mu}^{(k)} m_\mu
\,.
$$
Further, by the unitriangularity of $k$-Kostka numbers 
\eqref{trikostka}, we have a monomial expansion of 
dual $k$-Schur functions that can be seen as an alternative 
definition for these functions:

\begin{proposition} \label{propodualK}
For any $k$-bounded partition $\lambda$,
\begin{equation}
\mathfrak S_\lambda^{(k)} = m_\lambda + \sum_{\mu\lhd\lambda}
K_{\lambda\mu}^{(k)} m_\mu \,.
\end{equation}
\end{proposition}

This proposition reveals that the dual $k$-Schur functions are a basis
for the quotient of the symmetric function space by the ideal
$\mathfrak J^{(k)}$:

\begin{proposition}\label{propodualbasis}
The dual $k$-Schur functions form a basis of
$\Lambda/\mathfrak J^{(k)}$.
\end{proposition}

Recall that the $k$-Schur functions form a basis for 
$\Lambda/\langle h_i \, | \, i>k \rangle$.  The ideal 
$\mathfrak J^{(k)}$ is dual to $\langle h_i \, | \, i>k \rangle$ with 
respect to the scalar product defined on $\Lambda$ by 
$$
\langle h_\lambda,m_\mu\rangle = \delta_{\lambda\mu} \,.
$$
Since the definition of $k$-Schur function,
$s_{\lambda}^{(k)} = \sum_{\nu} \bar K_{\nu \lambda}^{(k)}\, h_{\nu}$,
implies that
\begin{equation*}
\langle s_\lambda^{(k)},\mathfrak S_\mu^{(k)} \rangle = 
\bigl \langle \, \sum_\alpha \bar K^{(k)}_{\alpha\lambda} h_\alpha,
\sum_\beta K_{\mu \beta}^{(k)} m_\beta
\bigr \rangle 
= \sum_\alpha   K_{\mu \alpha}^{(k)} \, \bar K^{(k)}_{\alpha\lambda} = 
\delta_{\lambda \mu} \, ,
\end{equation*}
as suggested by their name, the dual $k$-Schur 
basis is dual to the $k$-Schur basis.

\begin{proposition} 
\label{ortho}
Let $\lambda$ and $\mu$ be $k$-bounded partition.  Then,
$$
\langle s_\lambda^{(k)},\mathfrak S_\mu^{(k)}
\rangle = \delta_{\lambda\mu}
$$
\end{proposition}

\medskip

We can extract several combinatorial properties for dual $k$-Schur
functions from the $k$-Schur function properties using duality and
the following lemma.

\begin{lemma} \label{proposupermod}
 Let $f \in \Lambda^{k}$.  Then, for $g \in \Lambda$, we have
$$
\langle f\, , \, g \rangle = \langle f \, , \, g \mod \mathfrak J^{(k)} \rangle \, .
$$
\end{lemma}
\begin{proof}  
It suffices to consider $f=h_{\lambda}$, 
with $\lambda\in\mathcal P^k$.  If $A \in \mathfrak J^{(k)}$, 
then $A=\sum_\mu a_{\mu} \, m_{\mu}$ summing over
$\mu\not\in\mathcal P^k$. Thus, $\langle h_{\lambda}\, , \, A \rangle=0$
and the claim follows.
\end{proof}

\medskip

Since the $\omega$-involution is an isometry with respect to
$\langle \cdot,\cdot \rangle$, we discover from the 
action $\omega s_\mu^{(k)}= s_{\mu^{\omega_k}}^{(k)}$ \eqref{omegainv}
that $\omega$ acts naturally on the dual $k$-Schur functions.

\begin{proposition} \label{propoinvd}
Let $\lambda$ be a $k$-bounded partition.  Then
$$
\omega  \left( \mathfrak S_\lambda^{(k)} \right) \mod {\mathfrak J}^{(k)}
= \mathfrak S_{\lambda^{\omega_k}}^{(k)} \, .
$$
\end{proposition}

From the $k$-Pieri formula \eqref{kpieri},
with $f^\perp$ defined for any $f \in \Lambda$ by
$\langle  g, f^{\perp} \, h \rangle= \langle f\, g, h \rangle$,
we find
\begin{proposition}
\begin{equation}
h_\ell^{\perp} \, \mathfrak S_\nu^{(k)} = \sum_{\lambda\in \bar H_{\nu,\ell}^{(k)}}
\mathfrak S_\lambda^{(k)} \, ,
\end{equation}
where the sum is over partitions of the form:
$$
\bar H_{\nu,\ell}^{(k)}=
\left\{ \lambda \, \Big| \, 
\nu/\lambda={\rm horizontal~} \ell\text{-}{\rm strip} \quad {\rm and} \quad
\nu^{\omega_k}/\lambda^{\omega_k}={\rm vertical~} \ell\text{-}{\rm strip} 
\right\} \,.
$$
\end{proposition}

As with the Schur functions, the definition of $\mathfrak S_\lambda^{(k)}$ 
makes sense if $\lambda$ is replaced by a skew diagram.

\begin{definition} 
For $k$-bounded partitions $\mu\subseteq\nu$,
the ``dual skew $k$-Schur function" is defined by
\begin{equation}
\label{sumktab}
{\mathfrak S}_{\nu/\mu}^{(k)} = \sum_{\mathcal T} x^{k\text{-weight}
\, (\mathcal T)} \, ,
\end{equation}
where the sum is over all skew $k$-tableaux of shape $\core(\nu)/\core(\mu)$.
\end{definition}

The skew $k$-tableaux are well-defined since $\mu\subseteq\nu$ implies that 
$\core(\mu)\subseteq\core(\nu)$ (e.g.\cite{[LMcore]} Prop. 14).
The involution \eqref{tauinvo} on the set of $k$-tableaux, permuting the weight,
does the same to the set of skew $k$-tableaux.
Thus, $\mathfrak S_{\nu/\mu}^{(k)}$ is also a symmetric function
by the same reasoning that implies Proposition~\ref{dualsym}. Since
$K_{\nu/\mu,\lambda}^{(k)}$ denotes the number of skew $k$-tableaux 
of $k$-weight $\lambda$ and shape $\core(\nu)/\core(\lambda)$, 
we have the expansion:

\begin{equation}
\label{skewmon}
{\mathfrak S}_{\nu/\mu}^{(k)} = \sum_{\lambda} K_{\nu/\mu,\lambda}^{(k)} 
\, m_{\lambda} \, .
\end{equation}
Notice ${\mathfrak S}_{\nu/\mu}^{(k)} \in \Lambda/\mathfrak J^{(k)}$
since the sum is over $k$-bounded partitions $\lambda$ (a given
letter cannot have $k$-weight larger than $k$).  This form of the
skew dual $k$-Schur function makes it clear that the ``skew affine Schur 
functions" of \cite{[Lam]} are the same functions.

\medskip

Following the theory of usual skew Schur functions, we derive an
implicit formula for the skew dual functions in terms of dual
functions involving $k$-Littlewood-Richardson coefficients.

\begin{theorem}
For any $k$-bounded partitions $\mu\subseteq\nu$, 
$$
{\mathfrak S}_{\nu/\mu}^{{(k)}} = \sum_{\lambda} 
c_{\mu \lambda}^{\nu , k} {\mathfrak S}_\lambda^{(k)} \, .
$$
\end{theorem}
\begin{proof}
Since ${\mathfrak S}_{\nu/\mu}^{{(k)}}$ lies in
$\Lambda/\mathfrak J^{(k)}$, for which the
dual $k$-Schur functions form a basis,
$$
{\mathfrak S}_{\nu/\mu}^{{(k)}} = \sum_{\lambda} 
A_{\mu \lambda}^{\nu , k} {\mathfrak S}_\lambda^{(k)} \, ,
$$
for some $A_{\mu \lambda}^{\nu , k}$.
On one hand consider:
$$
\langle
s_\lambda^{(k)},\mathfrak S_{\nu/\mu}^{(k)}
\rangle
=
\langle
s_\lambda^{(k)},
 \sum_{\alpha} 
A_{\mu \alpha}^{\nu , k} {\mathfrak S}_\alpha^{(k)}
\rangle
=
A_{\mu \lambda}^{\nu , k}
$$
and
$$
\langle
s_\lambda^{(k)},\mathfrak S_{\nu/\mu}^{(k)}
\rangle
=
\langle
\sum_{\alpha} \bar K_{\alpha\lambda}^{(k)} h_\alpha,
\sum_{\beta} 
K_{\nu/\mu,\beta}^{(k)} m_\beta
\rangle
=
\sum_{\alpha} \bar K_{\alpha\lambda}^{(k)} 
K_{\nu/\mu,\alpha}^{(k)}
\,.
$$
On the other hand, since \eqref{eqhons} tells us
$h_{\lambda}\, s_{\mu}^{(k)} = \sum_{\nu} 
K_{\nu/\mu,\lambda}^{(k)} \, s_{\nu}^{(k)}$, we have
\begin{eqnarray*}
\langle
s_\mu^{(k)}s_\lambda^{(k)},\mathfrak S_{\nu}^{(k)}
\rangle
& = &
\langle
\sum_{\alpha} \bar K_{\alpha\lambda}^{(k)} h_\alpha s_\mu^{(k)},
\mathfrak S_{\nu}^{(k)}
\rangle
=
\langle
\sum_{\alpha} \bar K_{\alpha\lambda}^{(k)}
\sum_{\beta} K_{\beta/\mu,\alpha}^{(k)} s_\beta^{(k)},
\mathfrak S_{\nu}^{(k)}
\rangle
\\
& = &
\sum_{\alpha} \bar K_{\alpha\lambda}^{(k)}
K_{\nu/\mu,\alpha}^{(k)}
\end{eqnarray*}
and
$$
\langle
s_\mu^{(k)}s_\lambda^{(k)},\mathfrak S_{\nu}^{(k)}
\rangle
=
\langle
\sum_{\alpha} c_{\lambda\mu}^{\alpha, k} s_{\alpha}^{(k)},
\mathfrak S_{\nu}^{(k)}
\rangle = c_{\lambda\mu}^{\nu , k} 
\,.
$$
Therefore the result follows from
$$
A_{\lambda \mu}^{\nu,k} =\sum_{\alpha} \bar K_{\alpha\lambda}^{(k)}
K_{\nu/\mu,\alpha}^{(k)}=  c_{\lambda\mu}^{\nu , k} \, .
$$
\end{proof}

Given the duality between $\Lambda^k$ and $\Lambda/\mathfrak J^{(k)}$,
it is natural also to consider a skew $k$-Schur function.  The previous
proposition, exposing $k$-Littlewood-Richardson coefficients 
as the expansion coefficients for a the dual skew $k$-Schur function
in terms of dual $k$-Schur functions leads us to consider also
the coefficients in
$$
\mathfrak S_\lambda^{(k)} \, \mathfrak S_\mu^{(k)} = \sum_{\nu}
\mathfrak d_{\lambda \mu}^{\nu , k} \, \mathfrak S_\nu^{(k)} 
\mod \mathfrak J^{(k)} \, .
$$

\medskip

Similar to the $k$-Littlewood-Richardson coefficients,
Proposition~\ref{propoinvd} implies a symmetry satisfied by
these coefficients:
\begin{equation} \label{symd}
\mathfrak d_{\lambda \mu}^{\nu , k}  = \mathfrak 
d_{\lambda^{\omega_k} \mu^{\omega_k}}^{\nu^{\omega_k} , k}  \, .
\end{equation}
We also note that:
$$
\langle s_{\nu}^{(k)}, 
\mathfrak S_\lambda^{(k)} \, \mathfrak S_\mu^{(k)} \rangle
= \langle s_{\nu}^{(k)}, 
\mathfrak S_\lambda^{(k)} \, \mathfrak S_\mu^{(k)} \mod  \mathfrak J^{(k)}
\rangle = \mathfrak d_{\lambda \mu}^{\nu , k} \, .
$$

\medskip

We can now introduce the skew $k$-Schur function and discuss several
identities regarding the relations between these functions
and their dual.

\begin{definition}
For any $k$-bounded partitions $\mu\subseteq\nu$, the 
{\it ``skew $k$-Schur function"} is defined by
$$
s_{\nu/\mu}^{{(k)}} = \sum_{\lambda} 
\mathfrak d_{\mu \lambda}^{\nu , k} s_\lambda^{(k)} \, .
$$
\end{definition}

This given, our first property is:

\begin{proposition}  For any $f \in \Lambda$,
$$
\langle s_{\nu/\mu}^{(k)},f \rangle = 
\langle s_{\nu}^{(k)}, f \,  {\mathfrak S}_\mu^{(k)}  \rangle \, ,
$$
and for any $f \in \Lambda^{k}$,
$$
\langle f ,{\mathfrak S}_{\nu/\mu}^{(k)}\rangle = 
\langle f \, {s}_\mu^{(k)},{\mathfrak S}_{\nu}^{(k)}\rangle \, .
$$
\end{proposition}
\begin{proof}  From Proposition~\ref{propodualbasis} and 
Lemma~\ref{proposupermod},
it suffices to consider $f= \mathfrak S_\lambda^{(k)}$.
On one hand we have:
\begin{equation*}
\langle s_{\nu/\mu}^{(k)}, \mathfrak S_\lambda^{(k)}\rangle = 
\langle  \sum_{\alpha}
\mathfrak d_{\mu \alpha}^{\nu , k} s_\alpha^{(k)}, \mathfrak S_\lambda^{(k)} \rangle 
= \mathfrak d_{\mu \lambda}^{\nu , k} 
\,,
\end{equation*}
and on the other,
\begin{equation*}
\langle s_{\nu}^{(k)}, {\mathfrak S}_\lambda^{(k)} 
\, {\mathfrak S}_\mu^{(k)}\rangle  = 
\langle s_{\nu}^{(k)}, \sum \mathfrak 
d_{\mu \lambda}^{\alpha , k} {\mathfrak S}_\alpha^{(k)} \rangle 
=
\mathfrak d_{\mu \lambda}^{\nu, k} \, .
\end{equation*}
The second identity follows similarly.
\end{proof}

The $\omega$-involution again has a natural role in our study.
Given its action on $k$-Schur functions and their dual,
with the symmetries \eqref{symmetry} and \eqref{symd}, we find
\begin{proposition}
$$
\omega \left( {\mathfrak S}_{\nu/\mu}^{(k)} \right) \mod {\mathfrak J}^{(k)}= 
{\mathfrak S}_{\nu^{\omega_k}/\mu^{\omega_k}}^{(k)} \quad {\rm and} \quad
\omega \left( s_{\nu/\mu}^{(k)} \right) = s_{\nu^{\omega_k}/\mu^{\omega_k}}^{(k)} \, .
$$
\end{proposition}

The next proposition explains why the coproduct of $k$-Schur functions
in terms of $k$-Schur functions has the dual $k$-Littlewood-Richardson
coefficients as expansion coefficients, and thus connects the 
dual $k$-Schur functions with the cohomology of the loop Grassmannian
based on the conjecture of Shimozono.
Recall (e.g.  \cite{[Ma]}) that from the coproduct, 
$\Delta: \Lambda \to \Lambda(x) \otimes
\Lambda(y)$ by $\Delta f=f(x,y)$, a bialgebra structure is imposed:
$$
\langle \Delta f, g(x) \, h(y) \rangle = \langle f,gh \rangle \, ,
$$
where the first scalar product is in $\Lambda(x) \otimes
\Lambda(y)$.

\begin{proposition}
For any $\lambda\in\mathcal P^k$ and
two sets of indeterminants, $x$ and $y$,
$$
s_{\lambda}^{(k)}(x,y) = \sum_{\mu,\nu} \mathfrak d_{\mu \nu}^{\lambda,k}
s_{\mu}^{(k)}(x) \, s_{\nu}^{(k)}(y) =
 \sum_\nu
s_{\lambda/\nu}^{(k)}(x) \, s_{\nu}^{(k)}(y) \, ,
$$
and
$$
\mathfrak S_{\lambda}^{(k)}(x,y) = 
\sum_{\mu,\nu}  c_{\mu \nu}^{\lambda,k} \, 
\mathfrak S_{\mu}^{(k)}(x) \, \mathfrak S_{\nu}^{(k)}(y)
= \sum_\nu
\mathfrak S_{\lambda/\nu}^{(k)}(x) \, \mathfrak S_{\nu}^{(k)}(y) \, .
$$
\end{proposition}
\begin{proof}
For the first identity, given 
$\bigl \{ s_{\lambda}^{(k)} \bigr \}_{\lambda_1 \leq k}$  
forms a basis of $\Lambda^{k}$ and
$h_i(x,y) = \sum_{\ell=0}^i h_{i-\ell}(x) h_{\ell}(y)$,
we can assume
$$
s_{\lambda}^{(k)}(x,y) =
 \sum_{\gamma,\delta} e^{\lambda}_{\gamma\delta}
\, s_{\gamma}^{(k)}(x)\,s_{\delta}^{(k)}(y) \, ,
$$
for some $e^\lambda_{\gamma\delta}$.  We thus have
$$
\langle s_{\lambda}^{(k)}(x,y) , {\mathfrak S}_\nu^{(k)}(x) \, 
{\mathfrak S}_\mu^{(k)}(y) \, \rangle = 
\langle \sum_{\gamma,\delta} e^{\lambda}_{\gamma\delta}
s_{\gamma}^{(k)}(x)\,s_{\delta}^{(k)}(y)  , {\mathfrak S}_\nu^{(k)}(x) \, 
{\mathfrak S}_\mu^{(k)}(y) \, \rangle =e^{\lambda}_{\nu\mu} \, ,
$$
with the scalar product being taken in $\Lambda(x) \otimes
\Lambda(y)$.
On the other hand, using the definition of coproduct,
$$
\langle s_{\lambda}^{(k)}(x,y) , {\mathfrak S}_\nu^{(k)}(x) \, 
{\mathfrak S}_\mu^{(k)}(y) \, \rangle = \langle s_{\lambda}^{(k)}, 
{\mathfrak S}_\nu^{(k)} \, 
{\mathfrak S}_\mu^{(k)} \rangle =
 \langle s_{\lambda}^{(k)}, 
\sum_{\delta} \mathfrak d_{\nu \mu}^{\delta,  k}
{\mathfrak S}_\delta^{(k)} \rangle = 
\mathfrak d_{\nu \mu}^{\lambda , k} \, , 
$$
where we made use of Lemma~\ref{proposupermod} in the second equality.
Therefore,
$$
s_{\lambda}^{(k)}(x,y) =
 \sum_{\gamma,\nu} \mathfrak d^{\lambda\, (k)}_{\gamma\nu}
\, s_{\gamma}^{(k)}(x)\,s_{\nu}^{(k)}(y)  = \sum_\nu
s_{\lambda/\nu}^{(k)}(x) \, s_{\nu}^{(k)}(y) \, .
$$

For the second identity, we follow the same argument, using
the fact that if $\lambda_1 \leq k$, then
$m_{\lambda}(x,y)=\sum_{\mu,\nu \, : \, \mu_1,\nu_1\leq k } 
a_{\mu \nu} \, m_\mu(x) \, m_{\nu}(y)$.
\end{proof}

More generally, following the proof of
of (5.10) in \cite{[Ma]}, the preceding results imply:

\begin{proposition}
For any $\lambda\in\mathcal P^k$ and
two sets of indeterminants, $x$ and $y$,
$$
s_{\lambda/\mu}^{(k)}(x,y) = \sum_\nu
s_{\lambda/\nu}^{(k)}(x) \, s_{\nu/\mu}^{(k)}(y)
$$
and
$$
\mathfrak S_{\lambda/\mu}^{(k)}(x,y) = \sum_\nu
\mathfrak S_{\lambda/\nu}^{(k)}(x) \, \mathfrak S_{\nu/\mu}^{(k)}(y) \, .
$$
\end{proposition}

We conclude our exploration of the $k$-Schurs and their dual by
mentioning that they give rise to a refined Cauchy formula,
relying on the ideal $R^{(k)}$, generated by 
the indeterminants $y_i^{k+1}$ for $i=1,2,\dots$.  Note that moding out by $R^{(k)}$ 
simply amounts
to setting $y_i^{n}=0$ in the Cauchy kernel whenever $n>k$.

\begin{theorem}
Consider two bases of homogeneous symmetric functions,
$\left\{a_\lambda\right\}_{\lambda\in\mathcal P^k}$ 
and $\left\{b_\lambda\right\}_{\lambda\in\mathcal P^k}$ for
$\Lambda^{k}$ and $\Lambda/\mathfrak J^{(k)}$, respectively.
$$
\prod_{i,j} (1-x_i y_j)^{-1} \mod R^{(k)}= \sum_\lambda a_\lambda(x) b_\lambda(y) 
$$
iff $\langle a_\lambda , b_\mu \rangle=\delta_{\lambda \mu}$ for all $k$-bounded
partitions $\lambda$ and $\mu$.
\end{theorem}
\begin{proof}  The proof is similar to that of \cite{[Ma]}, I\,(4.6).
\end{proof}

\begin{corollary}
$$
\prod_{i,j}(1-x_i y_j)^{-1} \mod R^{(k)} =
\sum_{\lambda_1\leq k} h_\lambda(x) m_\lambda (y) 
=
\sum_{\lambda_1\leq k} s_\lambda^{(k)}(x) \mathfrak S_\lambda^{(k)}(y) \, .
$$
\end{corollary}

\section{Further work}

As detailed in the introduction, the Schur functions provide
a vehicle to directly reach the structure constants for 
multiplication in the cohomology of the Grassmannian from the
Littlewood-Richardson coefficients.  More generally, 
Theorem~\ref{gromovwitten} implies that the $k$-Schur functions
provide the analogous link between the quantum structure constants 
(or 3-point Gromov-Witten invariants) and the $k$-Littlewood-Richardson 
coefficients.  There are beautiful combinatorial methods for computing 
the Littlewood-Richardson coefficients that use, for example, skew 
tableaux or reduced words for permutations.  Although there has been 
progress in certain cases \cite{[BKT],[BMW],[BKMW],[KMSW],
[KTW],[SS1],[SS2],[Tu]}, 
a combinatorial interpretation for the 3-point Gromov-Witten invariants in 
complete generality remains an open problem.  The theory of $k$-Schur
functions suggests a number of natural approaches to this 
problem, with an extended notion of skew tableaux to $k$-skew tableaux, 
and the revelation \cite{[LMcore]} that affine permutations are the 
appropriate extended notion of permutations in this study.

\medskip

The conjecture that the (dual) $k$-Schur functions are the Schubert 
basis for the (cohomology) homology of the loop Grassmannian remains 
to be proven.  More generally, this supports the idea  that
there exists an affine version of Schubert polynomials related to $k$-Schur
functions.  In particular, the  (dual) $k$-Schur functions are indexed by $k$-bounded
partitions, which are in bijection with affine permutations in the quotient 
\cite{[BB]}.  Such affine permutations can be considered as the Grassmannian
version of affine permutations.  The results here, with the conjectures that 
the (dual) $k$-Schur basis is related to the Schubert basis for the 
(cohomology) homology of the loop Grassmannian,
suggest that these functions provide the the symmetric component of (dual) affine 
Schubert polynomials.  The first step in this direction is being
developed by Thomas Lam in a forthcoming paper \cite{[Lam]}
on affine Stanley symmetric functions.

\medskip

{\bf Acknowledgments}
{\it We thank Frank Sottile for his question and Mark Shimozono for 
his discussions and conjectures about the loop Grassmannian.  This 
work has greatly benefited from their interest.}

\end{document}